\documentstyle[12pt,amstex,psfig]{article}

\newtheorem{thm}{Theorem}
\newtheorem{lemma}[thm]{Lemma}
\newtheorem{cor}[thm]{Corollary}
\newcommand{\R}{{\Bbb R}}

\newcommand{\Z}{{\Bbb Z}}

\begin{document}
\title{A note on tiling with integer-sided rectangles.}
\author{Richard Kenyon\thanks{CNRS  UMR 128, Ecole Normale Sup\'erieure de Lyon,
46, all\'ee d'Italie, 69364 Lyon, France. Research at MSRI supported in part by
NSF grant no.DMS-9022140.}}
\date{}
\maketitle

\begin{abstract}
We show how to determine if a given simple rectilinear polygon can be tiled
with rectangles, each having an integer side.
\end{abstract}

\section{Introduction}
In \cite{Wag}, Stan Wagon provides us with 14 proofs of the fact that
if a rectangle $R$ is tiled with rectangles, each having at least one side
of integral length, then $R$ has a side of integral length.

Rather than simply add a fifteenth proof to his list, we would like
to address the much more general problem of tileability of arbitrary
rectilinear polygons.
In particular we will give an algorithm for deciding when a
rectilinear polygon (that is, a polygon with sides parallel to the axes)
can be tiled with rectangles, each having an integer side.

The proof idea comes from a neat method, due to John Conway in
the case of polyominos, for dealing
with tiling problems of this sort.
That method is to define a group, the ``tiling group''
which depends on the set of tiles involved in the particular problem,
and gives a necessary condition for the tileability of a given
simply connected region.

This method has been used with success in many cases where the tiles
have simple shapes (\cite{Thu}, \cite{Kpar}, \cite{KK}, \cite{Ksq}).
Here we apply it to Wagon's problem.
Surprisingly, the method resembles closely that of \cite{Thu} (or \cite{KK}) 
in the case of dominos.

One corollary to our construction is a simple description
of the set of all possible tilings of a polygon $R$ (Theorem \ref{rotate}).

\section{The tiling group}
Rather than define the tiling group we will simply use one of its quotient groups;
this will be sufficient for our purposes. For background into defining
the tiling group in a more general setting see the references.

Let $S^1$ denote the group $\R/\Z$.
Let $G=S^1*S^1$, the free (nonabelian) product of $S^1$ with itself.
Don't let the size or topology on this group worry you. Despite being
non-locally compact, it is a rather easy group to work with.

An element of $G$ is a product of elements in one or the other
factor, for example a typical element is
$$h(\frac13)h(\frac 12)v(-.12)h(\pi)v(\sqrt{2}),$$
where the symbols $h(\cdot), v(\cdot)$
refer to elements of the first or second free factor, respectively. 
An element of $G$ can be written in {\bf reduced form} using the 
(confluent) identities
$h(t_1)h(t_2)=h(t_1+t_2)$, $v(t_1)v(t_2)=v(t_1+t_2)$,
and $h(0)=v(0)=e$ (the identity). The reduced form of an element
is either $e$ or an alternating sequence of $h(t)$'s and $v(t)$'s, where the $t\not=0$
(the element can start and end in either an $h$ or a $v$).

Multiplication in $G$ is just concatenation of the corresponding expressions.

To a rectilinear polygonal path $\alpha$ starting from the origin
we associate an element of $G$ as follows. The path $\alpha$
is an ordered sequence of horizontal and vertical edges; to a horizontal
edge of ``length'' $t\in\R$ (that is an edge from a point with $x$-coordinate $a$
to a point with $x$-coordinate $t+a$) we associate the element $h(t)$.
For a vertical edge from $y$-coordinate $b$ to $y$-coordinate $b+t$
we associate the element $v(t)$. Then the element of $G$ corresponding to
a path is just the product of the elements coming from
the sequence of horizontal and vertical edges of that path.

For example, the path running counterclockwise 
around the boundary of a $t_1\times t_2$
rectangle with the origin in its lower left corner is:
\begin{equation}
h(t_1)v(t_2)h(-t_1)v(-t_2).\label{rect}
\end{equation}

For a simple rectilinear polygon $R$ having the origin on its boundary,
let $\alpha(R)$ be the
element of $G$ corresponding to the path running counterclockwise around the boundary
of $R$.

\begin{lemma} If $R$ is tileable with rectangles, each having an integer side,
then $\alpha(R)=e$ in $G$.
\label{trivialinG}
\end{lemma}

\begin{figure}[htb]
\centerline{\psfig{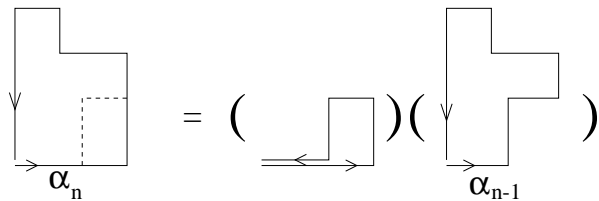}}
\caption{Inductive step in the proof that a tiling gives a product of 
lassos.\label{lassos}}
\end{figure}

\noindent{\bf Proof.} A tiling of $R$ gives a representation of $\alpha(R)$
as a product of ``lassos'', i.e. words of the form $xh(t)v(t')h(-t)v(-t')x^{-1}$
(conjugates of $[h(t),v(t')]$), where the words $[h(t),v(t')]$ are the boundaries
of the tiles: this can be proved inductively on the number of tiles, taking
at the inductive step a tile touching $\partial R$ (Figure \ref{lassos}),
see also \cite{CL}. 

Since each tile has an integral side, for
each lasso $x[h(t),v(t')]x^{-1}$, one of $t,t'$
is integral and so the lasso is trivial in $G$. 
Thus $\alpha(R)$ is the product of trivial elements of $G$ and so is trivial.
\hfill{$\Box$}
\medskip

This lemma is a necessary condition for tileability of a region by integer-sided
rectangles. It is not sufficient, as can be seen by considering for example
the (untileable) polygon whose boundary word is
$$h(\frac23)v(\frac13)h(-\frac13)v(1)h(-\frac23)v(-\frac13)h(\frac13)v(-1)=e.$$

\begin{cor} If a $t_1\times t_2$ rectangle $R$ is tileable with integer-sided
rectangles, then either $t_1\in\Z$ or $t_2\in\Z$.
\end{cor}

\noindent{\bf Proof.} The element $\alpha(R)=h(t_1)v(t_2)h(-t_1)v(-t_2)$
is trivial if and only if
either $h(t_1)=e$ or $v(t_2)=e$, that is, if and only if $t_1=0\bmod 1$
or $t_2=0\bmod 1$.
\hfill{$\Box$}
\medskip

We leave the reader to compare this proof with those of \cite{Wag}.

\section{The algorithm}
We give here an algorithm for deciding if a given simply rectilinear
polygon $R$ can be tiled with integer-sided rectangles.

Firstly, it is necessary that $\alpha(R)=e$ in $G$ by Lemma \ref{trivialinG}.

On the group $G$ define the distance function: for $x,y\in G$,
$d(x,y)=d(xy^{-1},e)$, where $d(x,e)$
is the length of the reduced word representing $x$ in terms of the number
of $h$'s and $v$'s. For example, the word in (\ref{rect}) has length $4$
assuming neither $t_1$ nor $t_2$ is integral.

Let $x_0\in G$ be the following base point: 
$$x_0=(h(-\frac12)v(-\frac12))^k=h(-\frac12)v(-\frac12)h(-\frac12)\ldots 
v(-\frac12),$$
where $k$ is sufficiently 
large so that $x_0$ is far away from $\alpha(R)$ and any of its prefixes.
We have $d(x_0,e)=2k$.

We define a height function $h_R$ on $\partial R$ as follows. This is a function on
points of $\partial R$ taking nonnegative integer values.
Let $y\in\partial R$, and
$\alpha_y$ be the path around $R$ counterclockwise
from $0$ to $y$.
Then define $h_R(y)=d(\alpha_y,x_0)$. 

Since $\alpha(R)=e$, we have $h_R(y)=d(\alpha_y',x_0)$ where $\alpha_y'$ is the path
from $0$ to $y$ going around $R$ clockwise instead ($\alpha_y$ and $\alpha_y'$ are the
same element of $G$).
Similarly,
if $R$ is tiled by integer-sided rectangles, we can extend $h_R$ to a height
function on all of the edges in the tiling: define $h_R$ on a point $y$ on an
interior edge
by simply taking any path $w_y$ along tile boundaries from the origin to that point;
then $h_R(y)\stackrel{def}{=}d(w_y,x_0)$ 
is independent of the path $w_y$ taken since for each tile
the path around its boundary is trivial in $G$.
\medskip

\noindent{\bf Remark.} Note that on an edge 
$z$ almost all the points have the same $h_R$-value.
The exceptions are points whose $h_R$-value is one less; these points
occur at distance exactly $1$ apart on $z$.

To tile the region $R$, we use the following result.
\begin{lemma} If $R$ is tileable, there is a tiling in which the maximum
of $h_R$ occurs on the boundary of $R$.
\label{maxonbdy}
\end{lemma}

\noindent{\bf Proof.} 
Suppose we have a tiling and the maximum is at a point $x$ on an
edge $z$ in the interior of $R$.
Suppose without loss that $z$ is vertical. Let $z'$ be a horizontal 
edge meeting $z$ at a point $y$. If $h_R(y)=h_R(x)$ then points on $z'$
are higher than $x$, contradicting our choice of $x$. So $h_R(y)=h_R(x)-1$.

By the remark which precedes this lemma,
every horizontal edge meeting $z$ is at an integer distance
from $y$. In particular the gap between adjacent horizontal edges
has integral length. So each rectangle adjacent to $z$ has integral height.

\begin{figure}[htb]
\centerline{\psfig{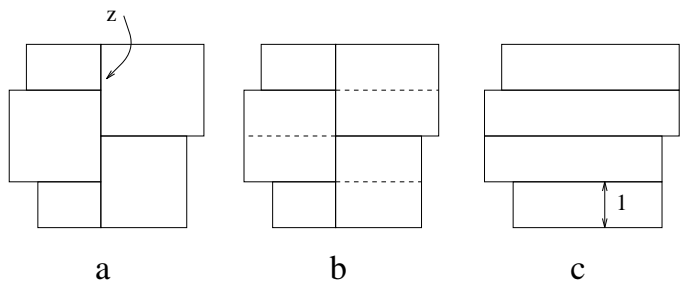}}
\caption{Removing an interior maximum.\label{maxedge}}
\end{figure}

Subdivide these rectangles using horizontal segments at each
height which is an integral distance from $y$ (Figure \ref{maxedge}b). 
This divides each
rectangle adjacent to $z$ into other rectangles of height $1$.
Now remove the edge $z$ completely; this has the effect of combining
pairs of rectangles left and right
of $z$ into single rectangles of height $1$ (Figure \ref{maxedge}c).

So we have a new tiling with integer-sided rectangles. On each edge we
added, the points are at a strictly lower height than $x$. 

By repeating this process for each interior edge which contains
a point of maximal height,
we eventually decrease the maximum in the interior. 
Since the maximum takes integer values,
eventually we find a tiling with the maximum on the boundary.
Furthermore there is a tiling in which points on interior edges 
are {\em strictly} lower than the boundary maximum.
\hfill{$\Box$}
\medskip

This lemma gives us an algorithm for tiling. 
First compute $h_R$ on the boundary of $R$. 
By Lemma \ref{trivialinG}, it is well-defined.
Find the maximum and minimum of $h_R$ on the boundary of $R$.
If the maximum minus the minimum is $1$, then all edges have integral length,
and so it is easy to tile. 

Suppose the maximal height difference is at least $2$.
By the remark preceding Lemma \ref{maxonbdy}, the maximum height 
occurs on a boundary edge $s$ of integer length. Assume without loss that $s$
is vertical. Put in as first tile a rectangle $R_1$ of side
$s$ and width $r$, where $r>0$ is the smaller of:
\begin{enumerate}
\item the smallest value for which
the height of all points on the edge $s'$ of $R_1$ opposite $s$ (extending $h_R$
to $R\cup R'$) is less than 
the maximum height, \item
the smallest value at which the edge $s'$ opposite $s$ touches the boundary of $R$ at a 
point on the interior of the edge $s'$. 
\end{enumerate}
We claim that if there is a tiling of $R$, there is a tiling
in which this rectangle $R'$ is a tile. By the lemma, there exists a tiling with the
maximum at $s$ whose interior edges are lower than $s$.
The rectangles adjacent to $s$ in this tiling must have integer height,
with their horizontal edges meeting $s$ at integer distance from the two
vertices of $s$ (since the horizontal edges must be lower than $s$). 
Furthermore these rectangles have widths at least as large as $r$ since their
vertical edges opposite $s$ either have height lower than $s$ or 
touch the boundary.
So in such a tiling we can cut off each of these tiles at distance $r$ from $s$
so that the added rectangle $R_1$ is a tile in a tiling of $R$.

So we can now continue with the remaining untiled
region $R'=R-R_1$. 
It may be that $R'$ has several components; however each component
is still simply connected, so
we run the algorithm in each component separately. 
On each component, for a tiling to exist the height function must be well-defined, that is,
the heights on new edges of $R_1$ must agree with the heights at points where they
touch $\partial R$. 

In either case $1$ or $2$ above in our choice of $r$, the total length
of boundary edges where the height is maximal has decreased by an integer amount:
in case $1$ it has decreased by the length of $s$, and in case $2$ it has decreased
by the length of $s'\cap\partial R$ which must be a positive integer.

So after a finite number of added rectangles the maximum on the
various boundaries has decreased. Since the minimum has not changed,
the maximal height difference
on each boundary has now decreased. Continue until this difference is equal
to $1$ on each component. Then the boundary must have
all edges of integral length and so can be trivially tiled. 

This completes the algorithm.
\label{alg}

\begin{figure}[htb]
\centerline{\psfig{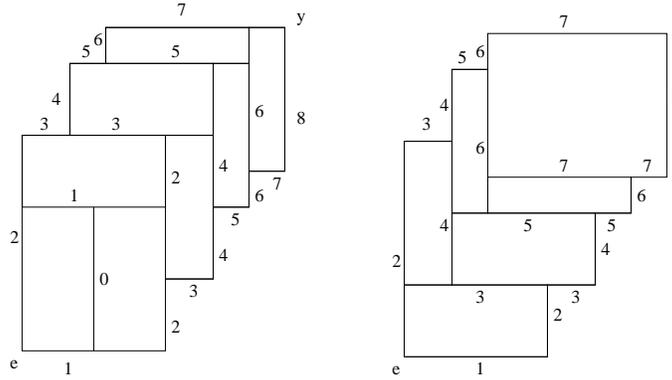}}
\caption{Lowest and highest tilings of a polygon.\label{poly}}
\end{figure}

\noindent{\bf Example.} 
In the polygon $R$ of Figure \ref{poly}, the counterclockwise path from $e$ to $y$
is the word 
$$h(1)v(\frac12)h(\frac13)v(\frac12)h(\frac14)v(\frac14)h(\frac14)
v(1),$$
and the clockwise path from $e$ to $y$ is
$$v(\frac32)h(\frac13)v(\frac12)h(\frac14)v(\frac14)h(\frac54),$$
these are equal in $G$ so $\alpha(R)=e$.

The heights (assuming the height of the first edge is $1$) are indicated
on the boundary (the height labelling an edge is the
largest height of points on that edge). 
The lowest tiling is shown in the first figure,
the highest in the second. The heights on interior edges are as shown.
Our algorithm would construct the first tiling.

A consequence of the local rearrangement of tiles as in Lemma \ref{maxonbdy}
and this algorithm is:
\begin{thm} Any two tilings of a polygon can be obtained from one another by
a sequence of the two types of operations: {\bf (1)} subdividing a tile
into two tiles along an edge of integer length, {\bf (2)} coalescing two tiles
which share a common edge of integer length into one tile by removing that edge.
\label{rotate}
\end{thm}

\noindent{\bf Proof.} 
Define a tiling $T_0$ which is that constructed by the algorithm, except 
don't stop when all the edges are integral but continue to apply the algorithm as before.
Once all the edges are integral or half-integral, they remain so by our
definition of base point $x_0$. So each rectangle added from this point on 
takes up at least area $1/2$, and thus the algorithm terminates after a finite
number of steps.

We claim any tiling $T$ can be transformed using operations (1) and (2)
to this tiling $T_0$.
By the lemma we 
can transform $T$ until the interior edges have heights strictly less than the boundary 
maximum.
Take a highest edge $s$ on $\partial R$ as before. 
We can now transform $T$ so that $R_1$, the rectangle added
in the algorithm, is a rectangle of $T$: each rectangle adjacent to $s$
has integral height, and so can be subdivided into two legal rectangles 
with a vertical edge at
distance $r$ from $s$. Then removing the horizontal edges which cut through $R_1$
leaves the tile $R_1$. 

The proof is then completed by induction on the number of tiles of $T_0$.
\hfill{$\Box$}

\section{Polyominos}
The preceding algorithm gives a bonus in the case of tiling a polyomino. 
Suppose we want to tile a simply connected 
polyomino with rectangles, each of which has a
side of length $n$.

The algorithm of section \ref{alg}
has the property that the other side-lengths of the tiles used
are in the lattice generated by the edge lengths of $R$. 
In the case of polyominoes, this
lattice is just $\Z$, and so the algorithm yields in this case a tiling with 
bars of length $1\times n$ and $n\times 1$. 

By stretching the vertical coordinate we can also similarly
solve the case of tiling with $1\times n$ and $m\times 1$ bars.
This case was previously solved by a similar algorithm in \cite{KK}.

Our algorithm runs in time linear in the area.

\section{Generalization}
Let $A$ be an arbitrary subgroup of $\R$ and let 
$S_A$ be the set of rectangles with at least one side in $A$.
Our algorithm also works to tile with tiles in $S_A$.
In this case the group $G=\R/A*\R/A$ loses all of its irrelevant topology
(which was implicit but not explicit in our understanding of $S^1*S^1$), but what
remains is the word metric coming from the distance function $d(x,y)$ which
is still well-defined.

The proofs of the lemmas extend easily to this case. For the rest of the algorithm,
if $A$ is not discrete then the choice of width $r$ of the added tile
is not well-defined. However we can take any $r$ sufficiently small so that the tile $R_1$
encounters no other vertex of $R$, and such that $s'$ has lower height than $s$.
We again continue with the untiled region until the boundary height difference
is $1$, in which case all the edges are in $A$, and so tiling is trivial.

\end{document}